\documentclass[11pt]{article}

\usepackage{euscript}

\evensidemargin\oddsidemargin

\usepackage{theorem,amsmath,amssymb}
\usepackage{graphics}
\usepackage{multicol}

\newtheorem{theo}{Theorem}[section]

\newtheorem{lemma}[theo]{Lemma}

\newcommand{\mysection}{
\section}

\def\a{\alpha}

\def\b{\beta}

\def\bhat{\widehat{\b}}
\def\Bhat{\widehat{B}}
\def\Bt{\widetilde{\beta}}

\def\D{\Delta}
\def\d{\delta}
\def\e{\varepsilon}

\def\E{\mathbb{E}}
\def\o{\omega}

\def\g{\gamma}

\def\l{\lambda}
\def\L{\mathcal{L}}

\def\Ib{\widehat{I}}
\def\N{\mathbb{N}}
\def\P{\mathbb{P}}
\def\p{\psi}

\def\R{\mathbb{R}}

\def\Rt{\widetilde{R}}

\def\s{\sigma}

\def\el{\overset{\L}{=}}

\def\EA{E_4}
\def\EB{E_5}
\def\EC{E_2}
\def\EG{E_{7}}
\def\EGb{E_{12}}
\def\EH{E_{13}} 
\def\EJ{E_{17}}
\def\EK{E_{10}}
\def\ELb{E_6}
\def\EM{E_{14}}
\def\EN{E_{15}}
\def\EO{E_{16}}
\def\EP{E_{19}}
\def\EQ{E_{20}}
\def\ER{E_{8}}
\def\ES{E_{18}}
\def\ET{E_{11}}
\def\EW{E_1}
\def\fn{E_3}
\def\EY{E_{9}}

\def\po{P_{\omega}}

\def\v1{v_{1}}

\def\bc{Borel--Cantelli }

\def\j1{J_{\Lambda}^<}

\def\k{\kappa}
\def\wk{W_{\k}}
\def\wkt{\widetilde{W}_{\k}}

\def\lB{L_B}

\def\lneg{L_X^{*-}(+\infty)}

\def\hbar{\overline{H}}

\def\rt{\widetilde{R}_{2+2\k}}
\def\rb{\widehat{R}_{2+2\k}}
\def\ty{\Lambda_Y}

\addtocounter{page}{-1}
\title{Almost sure asymptotics for a diffusion process in a drifted Brownian potential}

\author{Alexis Devulder
\footnote{Laboratoire de Probabilit\'es et Mod\`eles Al\'eatoires,
Universit\'e Paris VI, $4$ Place Jussieu, F-$75252$ Paris Cedex
$05$, France. E-mail: devulder@ccr.jussieu.fr.}}

\begin{document}

\maketitle

\begin{abstract}
We study a one-dimensional diffusion process in a drifted Brownian potential. We characterize the upper functions of its hitting times in the sense of Paul L\'evy, and determine the lower limits in terms of an iterated logarithm law.

\end{abstract}

\bigskip


\thispagestyle{empty}

\newpage


\mysection{Introduction}

We consider a diffusion process in random environment, defined as
follows. For $\k\in\R$, we introduce the random potential
\begin{equation}\label{eqdefinitionWk}
\wk(x):=W(x)-\k x/2,\qquad x\in\R,
\end{equation}

\noindent where $(W(x),\ x\in\R)$ is a standard two-sided
Brownian motion. We define a diffusion process ($X(t),\ t\geq
0)$ in the random potential $W_{\k}$, as solution to the
formal stochastic differential equation
$
\textnormal{d}X(t)=\textnormal{d}\b(t)-\frac{1}{2}\wk'(X(t))
\textnormal{d}t,
$ where $(\b(t),\ t\geq 0)$ is a Brownian motion independent of $W$
and $X(0)=0$. More rigorously, $X$ is a diffusion process such that
$X(0)=0$, and whose conditional generator given $W_{\k}$ is $
\frac{1}{2}e^{W_{\k}(x)}\frac{\text{d}}{\text{d} x}
\left(e^{-W_{\k}(x)}\frac{\text{d}}{\text{d} x}\right). $

We denote by $\po$ the law of $X$ conditionally on the environment
$\wk$, and call it the quenched law. We also define $\P(\cdot
):=\int \po(\cdot )P(\wk\in\text{d}\o)$, and call it the annealed
law.

The diffusion $X$, introduced by Schumacher (1985)
 and Brox (1986),
is generally considered as the continuous
time analogue of random walks in random environment (RWRE), which
have many applications in physics and biology.
For an account of general properties of RWRE, we refer to Zeitouni
(2004).

In this paper, we are interested in the transient case, that is,
we suppose $\k\neq 0$.
We may assume without loss of generality that $\k>0$. In this case,
$X(t)\to_{t\to+\infty}+\infty$ $\P$--a.s.

Our goal is to study the almost sure asymptotics of $X$.


We denote by $H$  the first hitting time of $r$ by $X$, that is,
\begin{equation}
    \label{H}
    H(r):=\inf\{t\geq 0,\quad X(t)>r\},\qquad r\geq 0.
\end{equation}

\noindent (See (\ref{HH}) for an analytic expression of $H(r)$). We
recall that there are three different regimes for $H$:

\medskip

\noindent {\bf Theorem A} {\it (Kawazu and Tanaka (1997))
 When $r$ tends to infinity,
 \begin{eqnarray}
     H(r)/r^{1/\k}
  & \overset{\L}{\longrightarrow}& c_0 S_{\k}^{ca} ,\qquad
     0<\k<1,\nonumber
     \\
     H(r)/(r\log r)
  & \overset{P.}{\longrightarrow}& 4, \qquad \quad\qquad \k=1,
     \nonumber
     \\
     H(r)/r
 & \overset{a.s.}{\longrightarrow}& 4/(\k-1), \qquad
     \k>1,\label{Kawazu_Tanaka_eq3}
 \end{eqnarray}
 where $c_0=c_0(\k)>0$ is a finite constant, the symbols
 $``\overset{\L}{\longrightarrow}"$,
 $``\overset{P.}{\longrightarrow}"$ and
 $``\overset{a.s.}{\longrightarrow}"$ denote respectively
 convergence in law, in probability and almost sure convergence,
 with respect to the annealed probability $\P$. Moreover,
 $S_\k^{ca}$ is a completely asymmetric stable variable of index
 $\k$, and is a positive variable for $0<\k<1$ (see
 \eqref{TransfoLaplace1} for its characteristic function).

}

\medskip

 In view of
 \eqref{Kawazu_Tanaka_eq3}, we only need to study the case $\k
\in (0,1]$. We prove

\begin{theo}
 \label{ThUpperLevyKinf1}
 Let $a(\cdot)$ be a positive nondecreasing function.
 If $0<\k<1$, then
 $$
  \sum_{n=1}^{\infty}\frac{1}{n a(n)} \
  \left\{
  \begin{array}{l}
      <+\infty
      \\
      =+\infty
  \end{array}
  \right.
  \Longleftrightarrow \limsup_{r\to\infty} \frac{H(r)}{[r
  a(r)]^{1/\k}}=
  \left\{
  \begin{array}{l}
      0
      \\
      +\infty
  \end{array}
  \right.
  \qquad\P\text{--a.s.}
 $$
 If $\k=1$, the statement holds under the additional
 assumption that $\limsup_{r\to+\infty} \frac{\log r}{a(r)}<\infty$.
\end{theo}

\begin{theo}
 \label{thLowerLevyKinf1}
 We have ($c_1(\k)\in (0,\infty)$ is given in equation \eqref{definitionc2}) :
 \begin{eqnarray}
    \liminf_{r\to+\infty} \frac{H(r)}{r^{1/\k} / (\log\log
    r)^{(1/\k)-1} }  & = &  c_1(\k) \qquad \P \text{--a.s.},\qquad
    0<\k<1,
    \label{eq_Th2_1}
 \\
    \liminf_{r\to+\infty} H(r)/(r\log r) & = & 4 \qquad
    \P \text{--a.s.},\qquad \k=1.
 \label{Th2_Keg1}
\end{eqnarray}
\end{theo}

\medskip

It was asked in Hu et al. (1999)
whether the
convergence in probability $H(r)/(r\log r) \to 4$ in Theorem A in the case $\k=1$ can be strengthened into
an almost sure convergence. Theorem \ref{ThUpperLevyKinf1} gives a
negative answer.

We observe that in the case $0<\k<1$, the process $H(\cdot)$ has the
same L\'evy classes as $\k$--stable subordinators (see Bertoin
(1996)
p.~92).

Theorems \ref{ThUpperLevyKinf1} and \ref{thLowerLevyKinf1} can be
stated for the process $X$ itself, by means of a standard
argument.

These results are useful for the study of the maximum local time of
$X$, see Devulder (2005).

\medskip

The rest of the paper is organized as follows.  We present in
Section~\ref{sectMainEstimates} some technical estimates which will
be needed later on; the proof of one of the technical estimates
(Lemma \ref{lemmaApproxLxLi}), is postponed until Section
\ref{sectannexe}. In Section \ref{SectEtudeHdeR}, we study the
L\'evy classes for the hitting times $H(r)$ and prove Theorems
\ref{ThUpperLevyKinf1} and \ref{thLowerLevyKinf1}.  Finally, we
prove Lemma \ref{lemmaApproxLxLi} in Section \ref{sectannexe}.

Throughout the paper, the letter $c$ with a subscript, denotes
unimportant constants that are finite and positive.


\mysection{Technical estimates}\label{sectMainEstimates}

We start by introducing $ A(x):=\int_0^x e^{\wk(y)}\text{d}y$ for
$x\in\R$,  which is a scaling function of $X$. We observe that,
since $\k>0$, $A(x){\rightarrow}_{x\to +\infty} A_\infty < \infty$.

For technical reasons, we have to introduce the random function
$F$ as follows. Fix $r>0$. Since the function $x\mapsto
A_{\infty}-A(x)$ is almost surely continuous and (strictly)
decreasing, there exists a unique
$F(r)\in\R$, depending only on the process $\wk$, such that
\begin{equation}
    A_{\infty}-A(F(r))=\exp(-\k r/2) =: \d(r).
    \label{e1p5}
\end{equation}

\noindent Our first technical estimate describes how close
$F(r)$ is to $r$, for large $r$.

\begin{lemma}
 \label{lemmaEncadrementFdeR}
 Let $\k>0$ and  $0<\d_0<1/2$. Define for $r>0$,
 \begin{equation}
    \label{eqEncadrementFdeR}
    \EW(r):=\{(1-5r^{-\d_0}/\k)r\leq
    F(r)\leq (1+5r^{-\d_0}/\k)r\}.
 \end{equation}
 Then for all large $r$,
 $
    \P(\EW(r)^c)\leq \exp(-r^{1-2\d_0}/4).
 $
 As a consequence, for any $\varepsilon>0$, we have, almost
 surely, for all large $r$,
 \begin{equation}
     \label{F}
     (1-\varepsilon) r \le F(r) \le (1+\varepsilon)r.
 \end{equation}
\end{lemma}

\medskip

\noindent{\bf Proof.} Let
$0<\d_0<1/2$, and fix $r>0$. Define $\wkt(u):=W(u+r)-W(r)-\k u/2$,
and $\widetilde{A}_{\infty}:=\int_0^{\infty}\exp(\wkt(u))\text{d}u$.
Hence, $\log[A_{\infty}-A(r)]=\log\widetilde{A}_{\infty}+\wk(r)$.
Let
$$
\EC(r):=\{(-2r^{-\d_0}-\k/2)r \leq \log[A_{\infty}-A(r)]\leq
(2r^{-\d_0}-\k/2)r\}.
$$

\noindent Recall that $\widetilde{A}_{\infty}\el 2/\g_\k$, where
$\g_\k$ is a gamma variable of parameter $\k$ (see e.g. Dufresnes,
(2000)), i.e., $\g_\k$ has density $e^{-x} x^{\k -1}/\Gamma(\k)$ for
positive $x$. Consequently,
\begin{eqnarray}
    \P(\EC(r)^c)
 \le \P[ \g_\k > 2e^{ r^{1-\d_0}} ] + \P [ \g_\k <2
    e^{-r^{1-\d_0} }] + \P[ |W(r)| > r^{1-\d_0} ]
 \le 3e^{ -r^{1-2\d_0}/2 }
    \label{eqLammaProbaInegDelta}
\end{eqnarray}

\noindent for $r$ large enough. Recall that
$A_{\infty}-A(F(r)) = \d(r) = e^{-\k r/2}$. On
$\EC[(1+5r^{-\d_0}/\k)r]$,
\begin{eqnarray*}
    \log \{A_{\infty}-A[(1+5r^{-\d_0}/\k)r]\}
 &\le&(\k r/2) ( -1 + 4 r^{-\d_0}/\k +
    o(r^{-\d_0}) ) (1+ {5r^{-\d_0} / \k})
    \\
 &<& \log [A_{\infty}-A(F(r))] ,
\end{eqnarray*}

\noindent where $f(r)=o(g(r))$ means $\lim_{r\to 0} {f(r)/ g(r)}
=0$. This gives the second inequality in \eqref{eqEncadrementFdeR}
by monotonicity of $A$.
Similarly, the first inequality holds on $\EC[(1-5r^{-\d_0}/\k)r]$.
This yields
$\P(\EW(r)^c)\leq \exp(-r^{1-2\d_0}/4)$ in view of
\eqref{eqLammaProbaInegDelta}.

Then, \eqref{F} follows from the Borel--Cantelli lemma and the
monotonicity of $F(\cdot)$.\hfill$\Box$

\medskip

With an abuse of notation, for $r\geq 0$, we denote by
$X\circ{\Theta_{H(r)}}$ the process $(X(H(r)+t)-r,\ t\geq 0)$,
which, conditionally on $\wk$, is a diffusion in the potential
$(\wk(x+r)-\wk(r),\ x\in\R)$, starting from $0$. Define
$H_{X\circ\Theta_{H(r)}}(s)=H(r+s)-H(r)$. Similarly,
$F_{X\circ\Theta_{H(r)}}$,
and $(H\circ F)_{X\circ\Theta_{H(r)}}$ denote respectively the
processes $F$
and $H\circ F$ for the diffusion $X\circ\Theta_{H(r)}$.
The following lemma is a modification of the \bc
lemma.

\begin{lemma}
 \label{lemmaBorelCantelli}
 Let $\k>0$.
 Let $(\D_n)_{n\ge 1}$ be a sequence of open sets in $\R$. Let
 $\a>0$, $r_n:=\exp(n^{\a})$ and $R_n:=\sum_{k=1}^n r_k$. If
$
     \sum_{n\geq 1} \P \left\{ (H\circ F)(r_{2n})
     \in \D_n \right\} = +\infty,
$
 then for any $\varepsilon>0$, almost surely, there exist
 infinitely many $n$ such that $H_{X \circ \Theta_{H(R_{2n-1})}} (t_n)
 \in\D_n$ for some
 $t_n \in [(1-\e)r_{2n}, (1+\e)r_{2n}]$.
\end{lemma}

\medskip

\noindent {\bf Proof.} Let $n\ge 1$, $x_n:=r_{2n-1}/2$, $\e_0>0$, $v_n:=2(\log n)/\k$ and
$$
\fn(n):=\left\{ \inf\nolimits_{t: \; H(R_{2n-1})\leq t\leq
H(R_{2n}+x_{n+1})} X(t)>R_{2n-2}+ x_n\right\}.
$$

\noindent
First, notice that
$
    \po(\fn(n)^c)
= ( 1+ [ \int\nolimits_{R_{2n-2} +x_n}^{R_{2n-1}}
    e^{\wk(x)} \text{d}x ]/[ \int_{R_{2n-1}}^{R_{2n}+x_{n+1}}
    e^{\wk(x)} \text{d}x] )^{-1}
$ a.s.
Define
\begin{eqnarray*}
    \EA
 := \{ \sup\nolimits_{0\leq x \leq r_{2n-1}-x_n } \left| \wk(
    x+ R_{2n-2}+x_n) - \wk(R_{2n-2} +x_n) +\k x/2 \right|
    \leq \e_0(r_{2n-1}-x_n)
    \},
    &\\
    \EB
 := \left\{ \sup\nolimits_{x\ge 0} [ \wk(x+R_{2n-1}) - \wk(R_{2n-1})] \leq v_n
 \right\}.\hphantom{AA.i...AAAAAAAAAAAAAAAAAAA}&&
\end{eqnarray*}
\noindent For large $n$, $
    \P(\EA^c)  \leq  2\exp[-\e_0^2(r_{2n-1}-x_n)/2]
$ and $
    \P(\EB^c)  =  \exp(-\k v_n)=1/n^2
$ (see Borodin et al. (2002), formula 1.1.4 (1)). Moreover, we have
for $n$ large enough, on $\EA\cap\EB$,
\begin{eqnarray}
    \po(\fn(n)^c)
 &\le& \k\frac{(r_{2n} +x_{n+1} ) \exp[ v_n+ \wk(R_{2n-1}
    )]} {\exp[ \wk( R_{2n-2} +x_n) -\e_0(r_{2n-1}
    -x_n)]}
    \nonumber
    \\
 &\le& \k( r_{2n} +x_{n+1} ) \exp[ v_n + (2\e_0-\k/2)
    (r_{2n-1} -x_n)].
\nonumber
\end{eqnarray}
Integrating this over $\EA\cap \EB$ yields
$\label{eqlemmaprobaFncBis} \sum_{n=1}^{+\infty}\P(\fn(n)^c)<\infty$
for $\e_0<\k/4$.

To complete the proof of Lemma \ref{lemmaBorelCantelli}, we define
\begin{eqnarray*}
    \mathcal{D}_n
 &:=& \{\exists t_n \in [(1-\e)r_{2n}, (1+\e)r_{2n}],
      H_{X\circ \Theta_{H(R_{2n-1})}} (t_n)
    \in \D_n
    \},\\
    \mathcal{E}_n
 &:=&\{(1-5r_{2n}^{-\d_0}/\k)r_{2n}\leq
    F_{X\circ\Theta_{H(R_{2n-1})}}
    (r_{2n})\leq(1+5r_{2n}^{-\d_0}/\k)r_{2n}\}.
\end{eqnarray*}

\noindent Let ${\widetilde t_n} := F_{X\circ\Theta_{H(R_{2n-1})}}
(r_{2n})$. We have $ \mathcal{D}_n \cap \fn(n) \supset \{ H_{X \circ
\Theta_{H(R_{2n-1})}} ({\widetilde t_n})  \in \D_n\} \cap
\fn(n)\cap\mathcal{E}_n . $ By assumption, $\sum_n \P \{
²H_{X\circ\Theta_{H(R_{2n-1})}}({\widetilde t_n}) \in\D_n \}
=\infty$. Moreover, $\P(\mathcal{E}_n)=\P(\EW(r_{2n}))$. Since
$\sum_{n=1}^{+\infty}\P(\fn(n)^c)<\infty$,  this  and Lemma
\ref{lemmaEncadrementFdeR}, yield $\sum_{n\in\N}\P(
\mathcal{D}_n\cap\fn(n))=+\infty$.

Since $\mathcal{D}_n\cap \fn(n)$, $n\geq 1$, are independent events,
the \bc lemma yields Lemma \ref{lemmaBorelCantelli}.~$\Box$

\medskip

In the rest of the paper, if $(\b(s),\ s\geq 0 )$ is a
Brownian motion, we denote its local time  by $(L_\b(t,x)\ t\geq
0,x\in\R)$, and define $\tau_\b(x):=\inf\{t>0,\ L_\b(t,0)=x\}$,
$x\geq 0$. For $v>0$, we define the Brownian motion $(\b_v(s),\
s\geq 0)$ by $\beta_v(s):=(1/v) \beta(v^2 s)$. We also introduce for
$\d_1>0$,
\begin{eqnarray}
    \l
 &:=& 4(1+\k) , \quad c_2:= 2(\lambda/\k)^{\d_1},
    \quad \psi_\pm (r):= 1\pm {c_2/ r^{\d_1}},
    \quad t_\pm (r):= {\k\psi_\pm(r)r / \l} ,
    \label{e1p22}
    \\
    K_\b({\k})
 & := & \int_0^{+\infty}x^{1/\k-2}L_{\b}(\tau_\b(\l),x)\text{d}x,
    \qquad 0<\k<1,
    \label{e3p18}
    \\
    C_\b
 & := & \int_0^1\frac{L_\b(\tau_\b(8),x)-8}{x}\text{d}x+
\int_1^{+\infty}\frac{L_\b(\tau_\b(8),x)}{x}\text{d}x.
    \label{e2p18}
\end{eqnarray}

\noindent We have the following approximation result.

\begin{lemma}
 \label{lemmaApproxLxLi}
 Let $0<\k\leq 1$ and  $\e\in (0,1)$. For $\d_1>0$ small enough, there exist
 $c_3>0$ and a
 Brownian motion $(\b(t),\ t\ge 0)$
 such that  for some $\a>0$ and
 all large $r$,
 $
 \P \{ \ELb(r) \} \ge 1- r^{-\alpha},
 $
 where
 \begin{eqnarray}
     \label{e4p25}
     \ELb(r)
  &:=& \{ (1-\e) \Ib_-(r)\le H(F(r))\le (1+\e)
     \Ib_+(r)\} ,
     \\
     \label{eqI''CasKinf1}
        \Ib_\pm(r)
  &:=&
     \left\{
     \begin{array}{lr} 4 \k^{1/\k-2}t_\pm(r)^{1/\k}
       \{ K_{\b_{t_\pm (r)}}(\k) \pm
       c_3t_\pm (r)^{1-1/\k}\},
      &  \quad 0<\k<1,
      \\[2mm]
       4t_\pm(r) \{C_{\b_{t_\pm (r)}} + 8 \log t_\pm
       (r) \},
      & \quad \k=1 .
     \end{array}\right.
 \end{eqnarray}
\end{lemma}

\medskip

\noindent {\bf Proof.} Postponed to Section
\ref{sectannexe}.\hfill$\Box$



\mysection{Proof of Theorems \ref{ThUpperLevyKinf1} and
\ref{thLowerLevyKinf1} }\label{SectEtudeHdeR}

In this section, we assume $0<\k\leq 1$, and prove Theorems
\ref{ThUpperLevyKinf1} and \ref{thLowerLevyKinf1}.


Let $S_\k^{ca}$ be a (positive) completely asymmetric stable
variable of index $\k$, and $C_8^{ca}$ a completely asymmetric
Cauchy variable of parameter $8$. Their characteristic functions are:
\begin{eqnarray}
    \E\exp(itS_\k^{ca})
 & = & e^{-|t|^\k\left(1-i\, \textnormal{sgn}(t)
    \tan(\frac{\pi\k}{2})\right) },
    \label{TransfoLaplace1}\qquad
    \E\exp(itC_8^{ca})
  =  e^{- 8 \left(|t|+i t \frac{2}{\pi} \log|t| \right)}.
\end{eqnarray}

Recall $\Ib_\pm$ from \eqref{eqI''CasKinf1}. By Biane and Yor
(1987), for $\l>0$,
\begin{eqnarray}
    \Ib_\pm(r)
    & \el & t_\pm(r)^{1/\k}\{ c_4 \, S_\k^{ca} \pm c_5 \,
    t_\pm(r)^{1-1/\k}\},\qquad\quad 0<\k<1,
    \label{eqDefItierceKinf1}\\
     \label{eqDefItierceKeg1}
    \Ib_\pm(r)
    & \el & 4t_\pm(r) [8c_6 +(\pi/2)
    C_8^{ca} + 8\log t_\pm(r) ],\qquad  \k=1,
    \\
    \psi(\k)
 &:=& \left(\frac{\pi \k}{4\Gamma^2(\k)\sin(\pi
    \k/2)}\right)^{1/\k},\qquad c_4 := 8
    \psi(\k)\l^{1/\k}\k^{-1/\k} ,
    \label{DefPsi}
\end{eqnarray}

\noindent and $c_5 >0$ and $c_6>0$ are unimportant constants.



\medskip

\noindent {\bf Proof of Theorem \ref{ThUpperLevyKinf1}.} 
Let $r_n:=e^n$ and $R_n:=\sum_{k=1}^n r_k$. Let $a(\cdot)$ be a
positive nondecreasing function. Without loss of generality, we
assume that $a(n) \to_{n\to+\infty} +\infty$.

We start with the case $0<\k<1$. By Samorodnitsky and Taqqu (1994)
p.~16), $\P(S_\k^{ca}>x)\sim_{x\to +\infty} c_7x^{-\k}$, where
$f(x)\sim_{x\to+\infty} g(x)$ means $\lim_{x\to+\infty}
f(x)/g(x)=1$, and $c_7$ is a positive constant depending on $\k$.

Recall $t_\pm(\cdot)$ from \eqref{e1p22}. By Lemma
\ref{lemmaApproxLxLi} and \eqref{eqDefItierceKinf1}, for large
$r$, we have
\begin{eqnarray}
    \P[H(F(r))>(a(e^{-2}r) t_+(r))^{1/\k}]
 &\le& c_8/a(e^{-2}r ) + (\log r )^{-2}.
    \label{e1p42}
\end{eqnarray}

Assume $\sum_{n\geq 1} {1\over a(r_n)}<\infty$. By the
Borel--Cantelli lemma, almost surely for $n$ large enough, $
H(F(r_n))\leq [a(r_{n-2})t_+(r_n)]^{1/\k} $. On the other hand, by
Lemma \ref{lemmaEncadrementFdeR}, almost surely for all large $n$,
we have $r_{n+1}\le F(r_{n+2})$, which implies that for $r\in [r_n,
r_{n+1}]$,
$$
H(r)\leq H[F(r_{n+2})] \leq [a(r_{n})t_+(r_{n+2})]^{1/\k} \leq
[\p_+(r_{n+2})\k r_{n+2}a(r_{n})/\l]^{1/\k} \leq c_9 [
ra(r)]^{1/\k}.
$$

\noindent Therefore, $ \limsup_{r\to+\infty}\frac{H(r)}{[r
a(r)]^{1/\k}}\leq c_9$ $\P\text{--a.s.},$ implying the ``zero'' part
of Theorem \ref{ThUpperLevyKinf1}, since we can replace $a(\cdot)$
by any constant multiple of $a(\cdot)$.


To prove the ``infinity" part, we assume $\sum_{n\geq 1} {1\over n
a(n)} = +\infty$, and observe that, by a similar argument leading
to \eqref{e1p42}, we have, for $r$ large enough,
\begin{equation}\label{displayP15a}
\P[H(F(r))>(a(e^2r)t_-(r))^{1/\k}]\geq c_{10}/a(e^2r)-(\log r)^{-2}.
\end{equation}

\noindent It follows from Lemma \ref{lemmaBorelCantelli} that $
\sup_{t\in[(1-\e)r_{2n},(1+\e)r_{2n}]}
H_{X\circ{\Theta_{H(R_{2n-1})}}}(t) >
[a(r_{2n+2})t_-(r_{2n})]^{1/\k}$ almost surely for infinitely many
$n$, which implies, for these $n$,
\begin{equation}
\sup\nolimits_{t\in[(1-\e)r_{2n},(1+\e)r_{2n}]}
H(R_{2n-1}+t)/[a(R_{2n-1}+t)(R_{2n-1}+t)]^{1/\k}\geq c_{11}
.\label{displayP15b}
\end{equation}

\noindent This gives $\limsup_{r\to+\infty}\frac{H(r)}{[r
a(r)]^{1/\k}}\geq c_{11}$ $\P$--a.s., which proves Theorem
\ref{ThUpperLevyKinf1} in the case $0< \k<1$.


It remains to treat the case $\k=1$. We recall that there exists a
constant $c_{12}>0$ such that $\P(C_8^{ca}>x){\sim}_{x\to +\infty}
c_{12}/x$ (see e.g. Samorodnitsky et al. (1994) p.~16). Hence,
\begin{eqnarray}
  \P \left\{ H(F(r))>4t_+(r)(1+\e) [8c_6 +a(e^{-2}r) +
    8\log t_+ (r) ] \right\}
 \le c_{12}\pi/ a(e^{-2}r) + (\log
    r)^{-2}
    \label{littlehope}
\end{eqnarray}

\noindent by Lemma \ref{lemmaApproxLxLi} and
\eqref{eqDefItierceKeg1}, for large $r$. Assume $\sum_{n\geq 1}
{1\over na(n)} <\infty$. Then by the \bc lemma, almost surely, for
all large $n$, $ H(F(r_n)) \le 4 t_+(r_n)(1+\e)  [8c_6 +a(r_{n-2})+
8 \log (t_+(r_n)].$ Under the additional assumption $\limsup_{r\to
+\infty}(\log r)/a(r)<\infty$, we have, almost surely, for all large
$n$ and $r\in [r_n, r_{n+1}]$ (thus $r \le F(r_{n+2})$ by Lemma
\ref{lemmaEncadrementFdeR}),
$$
H(r)\leq H(F(r_{n+2})) \le c_{13} r_{n+2} [ a(r_{n})+ \log r_{n+2}]
\le c_{14}r a(r) ,
$$

\noindent since $t_+(r)=\psi_+(r)\k r/\l$.  This yields the ``zero''
part of Theorem \ref{ThUpperLevyKinf1} in the case $\k=1$.


For the ``infinity" part, we assume $\sum_{n\geq 1} {1\over
na(n)} =+\infty$. As in \eqref{littlehope}, we have, for large
$r$,
$$
\P \left\{ H(F(r))>4t_-(r)(1-\e) a(e^2r) \right\} \geq \pi c_{12}/[4
a(e^2r)] -(\log r)^{-2}.
$$
As in the displays between \eqref{displayP15a} and
\eqref{displayP15b}, this yields the ``infinity" part of Theorem
\ref{ThUpperLevyKinf1} in the case $\k=1$. \hfill$\Box$


\medskip

\noindent {\bf Proof of Theorem \ref{thLowerLevyKinf1}.} We start with the case $0<\k<1$.
Recall (see Bertoin (1996)
\begin{equation}
    \label{bertoin}
    \log\P( S_\k^{ca} <x) {\sim}_{x\to 0,\ x>0} -
    c_{15} x^{-\k/(1-\k)},
\end{equation}

\noindent where $c_{15}$ is a constant depending only on $\k$.
Consequently, for $r$ large enough, by \eqref{eqDefItierceKinf1} and
Lemma \ref{lemmaApproxLxLi}, for any (strictly) positive function
$f$,
\begin{equation}
    \P[H(F(r)) < t_-(r)^{1/\k}f(r)]
 \le \exp \left[ -(c_{15}- \e) \left(
    \frac{ (1-\e) c_4 }{f(r)+ c_{16} r^{1-1/\k} }
    \right)^{\k/(1-\k)} \right] + \log^{-2} r.
    \label{eqH+inf1}
\end{equation}

\noindent Let $s_n:=\exp(n^{1-\e})$ and $ f(r):=\left(
\frac{1-\e}{1+\e}\right)^{(1-\k)/\k} \frac{(1-\e) c_4
(c_{15}-\e)^{(1-\k)/\k}}{(\log\log r)^{(1-\k)/\k}}- {c_{16} \over
r^{1/\k -1}}$.  As a consequence,  $\sum_n \P[ H(F(s_n)) <
t_-(s_n)^{1/\k} f(s_n)] < \infty$, which, by means of the \bc lemma,
implies that, almost surely, for all large $n$, $ H(F(s_n))\geq
t_-(s_n)^{1/\k}f(s_n)$.

Recall from Lemma \ref{lemmaEncadrementFdeR} that,
almost surely, for all large $n$, we have $F(s_n) \le
(1+\varepsilon) s_n$. Let $r$ be large. There exists $n$ (large)
such that $(1+\varepsilon) s_n \le r\le (1+2\varepsilon)s_n$. Then
$$
H(r) \ge H(F(s_n)) \ge  t_-(s_n)^{1/\k} f(s_n) \ge t_-^{1/\k}
\left[r/(1+2\e)\right] f[r/(1+2\e)] .
$$

\noindent Plugging the value of $t_-( {r\over 1+2\e})$ (defined in
\eqref{e1p22}), this yields inequality ``$\geq$" of \eqref{eq_Th2_1}
with
\begin{equation}
    \label{definitionc2}
    c_1(\k) = 8\psi(\k) c_{15}^{(1/\k)-1}
\end{equation}

\noindent where $c_{15}= c_{15}(\k)$ is defined in \eqref{bertoin},
and $\psi$ in \eqref{DefPsi}.


To prove the upper bound, let $
g(r):=\e+\left(\frac{1+\e}{1-\e}\right)^{(1-\k)/\k} \frac{(1+\e) c_4
(c_{15}+\e)^{(1-\k)/\k}}{(\log\log r)^{(1-\k)/\k}} + {c_{17} \over
r^{1/\k -1}}$, $r_n:=\exp(n^{1+\e})$ and $R_n:=\sum_{k=1}^n r_k$. By
means of an argument similar to the one leading to \eqref{eqH+inf1},
and Lemma \ref{lemmaBorelCantelli}, there exist almost surely
infinitely many $n$ such that
$$
\inf\nolimits_{u\in[(1-\e)r_{2n},(1+\e)r_{2n}]}
H_{X\circ{\Theta_{H(R_{2n-1})}}}(u)<[t_+(r_{2n})]^{1/\k}g(r_{2n}).
$$

\noindent In addition, by Theorem \ref{ThUpperLevyKinf1}, $
H(R_{2n-1})<\left[ R_{2n-1}\log^2 R_{2n-1}
\right]^{1/\k}\leq\e[t_+(r_{2n})]^{1/\k}g(r_{2n})$, almost surely,
for all large $n$,  since $R_k\le k \exp(-k^\e)r_{k+1}$ for all
large $k$, which yields
$$
\inf_{v\in[R_{2n-1}+(1-\e)r_{2n}, \, R_{2n-1}+(1+\e)r_{2n}]}
H(v)<(1+\e)[t_+(r_{2n})]^{1/\k}g(r_{2n}) .
$$

\noindent This gives inequality ``$\leq$" of \eqref{eq_Th2_1}, and
thus yields Theorem \ref{thLowerLevyKinf1} in the case $0<\k<1$.


We now assume $\k=1$ (thus $\l=8$). By Samorodnitsky and Taqqu
(1994) Proposition 1.2.12), $\E[\exp(-C_8^{ca})]=1$ (in the notation
of Samorodnitsky and Taqqu (1994) $S_1(8,1,0)$). Hence,
\begin{equation}
    \label{ProbaSamoTaqqu}
    \P[ C_8^{ca} \le -\varepsilon\log r ] \le r^{-\varepsilon}\,
    \E[ \exp( -C_8^{ca}) ] = r^{-\varepsilon}.
\end{equation}

\noindent By Lemma \ref{lemmaApproxLxLi} and
\eqref{eqDefItierceKeg1}, we have, for all large $r$,
\begin{equation*}
  \P \left\{H(F(r))\leq 32 t_-(r) (1-2\varepsilon) [ c_6
    + \log t_-(r)] \right\}
 \le \P\big( C_8^{ca}\le -
    \frac{16\e\log t_-(r)}{\pi(1-\e)}
    \big)
    +{1\over\log^{2}r}
    \le r^{-c_{18}}.
\end{equation*}

\noindent Let $s_n := \exp(n^{1-\varepsilon})$. Thus, by the \bc
lemma, almost surely, for all large $n$, $ H(F(s_n)) > 32 t_-(s_n)
(1-2\varepsilon) [c_6 + \log t_-(s_n)]$, which is greater than
$4(1-3\e)s_n\log s_n$. In view of (the last part of) Lemma
\ref{lemmaEncadrementFdeR}, this yields inequality ``$\geq$" in
\eqref{Th2_Keg1}. The inequality ``$\le$'', on the other hand,
follows immediately from Theorem A (that $H(r)/(r\log r)\to 4$ in
probability). Theorem \ref{thLowerLevyKinf1} is proved.\hfill$\Box$



\mysection{Proof of Lemma \ref{lemmaApproxLxLi}}
 \label{sectannexe}

Let $\k>0$. Recall $A(x)=\int_0^x e^{\wk(u)} \text{d}u$, and
$A_{\infty}=\lim_{x\to+\infty} A(x)<\infty$, $\P$--a.s. For any
Brownian motion $\b$ and any $r>0$, we define $\s_\b(r):=\inf\{t>0,\
\b(t)=r\}$.

Following Hu et al. (1999), there exists a Brownian motion $B$
independent of $W$ such that
\begin{equation}
    H(r)
  =  \int_{-\infty}^0 + \int_0^{A(r)} e^{-2\wk [A^{-1}(x)]}
    L_B\{\s_B[A(r)], x\} \text{d}x
    :=H_-(r)+H_+(r),
    \label{HH}
\end{equation}

Recall $F$ from \eqref{e1p5} and notice that $F(r)>0$ on $\EW(r)$.
Let $\Delta(r):=\d(r)^{-1}A(F(r))$. Following Hu et al. (1999)
p.3930), there exists two Bessel processes $R_2$ and $R_{2+2\k}$, of
dimensions $2$ and $(2+2\k)$ respectively, starting from $0$, such
that
$$
    H_+(F(r))
=
\int_0^{A(F(r))}\frac{16R_2^2(v)\textnormal{d}v}{R_{2+2\k}^4(v+\d(r))}
=
\int_0^{\Delta(r)}\frac{16R_2^2(\d(r)v)\textnormal{d}v}{R_{2+2\k}^4[\d(r)(1+v)]}
=
\int_0^{\Delta(r)}\frac{16\widetilde{R}_2^2(v)\textnormal{d}v}{\rt^4(1+v)},
$$
where $R_2=^{\L}\widetilde{R}_2$ and
$R_{2+2\k}=^{\L}\widetilde{R}_{2+2\k}$. As in Hu et al (1999), we
define  a Jacobi process of dimension $(d_1,d_2)$ as the solution of
\begin{equation}
    \label{eqDefJacobi}
    \text{d} Y(t) = 2\sqrt{Y(t) (1-Y(t))} \text{d} \bhat(t) +
    [d_1-(d_1+d_2)Y(t)] \text{d} t,
\end{equation}

\noindent where $\bhat$ is a standard Brownian motion. According to
Warren and Yor (1997), there exists a Jacobi process $(Y(t),\ t\geq
0)$ of dimension $(2,2+2\k)$, starting from $0$, such that for any
$u\ge 0$,
\begin{equation}
\frac{\widetilde{R}_2^2(u)}{\widetilde{R}_2^2(u)+\rt^2(u+1)} =
Y\circ \ty(u), \qquad
    \label{defThetaY}
    \ty(u):=\int_0^u\frac{\text{d}s}{\widetilde{R}_2^2(s)+\rt^2(s+1)}.
\end{equation}

\noindent In particular, $(\ty(t),\ t\geq 0)$ is independent of
$Y$. As a consequence, for all $r\geq 0$,
\begin{eqnarray}
H_+[F(r)] & = & 16\int_0^{\d(r)^{-1}A(F(r))} \frac{[Y\circ\ty(u)]
    \ty'(u) \text{d}u}{[1-Y\circ \ty(u)]^2}
    =16\int_0^{\g(r)} \frac{Y(u)}{(1-Y(u))^2} \text{d}u
\nonumber\\[-1mm]
    \label{e10p7}
    \g(r) & := & \ty[\d(r)^{-1}A(F(r))].
\end{eqnarray}

Let $\a_{\k}:= 1/(4+2\k)$ and let $T_Y(\a_\k):=\inf\{t>0,\
Y(t)=\a_{\k}\}$. Define
\begin{equation}
    \label{e1p7}
    \begin{array}{ll}
    \qquad \hbar(r):= 16\int\limits_0^{T_Y(\a_\k)}
    \cfrac{Y(u)}{(1-Y(u))^2}\text{d}u,
   & \qquad
    H_0(r):=16\int\limits_{T_Y(\a_\k)}^{\g(r)}
    \cfrac{Y(u)}{(1-Y(u))^2}\text{d}u,
    \end{array}
\end{equation}

\noindent and notice that
\begin{equation}
    H_+(F(r))=\hbar(r)+H_0(r),
   \qquad
    \{T_Y(\a_\k)\leq 64\log r\}\subset
    \{\hbar(r)\leq c_{19}\log r\}.
    \label{eqMajorAvantAlphaK}
\end{equation}

\noindent Observe that a scale function of $Y$ is $S(y):=
\int_{\a_\k}^y \frac{\text{d}x}{x(1-x)^{1+\k}}$. There exists a
Brownian motion $(\b(t),\ t\geq 0)$ such that for all $t\geq 0$,
$$
Y[t+T_Y(\a_\k)]=S^{-1}\{\b[U(t)]\}, \qquad  U(t):=4 \int_0^t
\frac{\text{d}s}{Y[s+T_Y(\a_\k)] \{1-Y[s+T_Y(\a_\k)] \}^{1+2\k}}.
$$

The rest of the proof of Lemma \ref{lemmaApproxLxLi} requires some
more estimates, stated as Lemmas
\ref{lemmaApproxTy}--\ref{lemmaApproxJCasKinf1} below. We admit
these lemmas for the moment, and complete the proof of
Lemma~\ref{lemmaApproxLxLi}.

\begin{lemma}\label{lemmaApproxTy}

Let $(R(t),\ t\geq 0)$ be a Bessel process of dimension $d>4$,
starting from $R_0\el\widetilde{R}_{d-2}(1)$, where
$(\widetilde{R}_{d-2}(t),\ t\in[0,1])$ is a $(d-2)$--dimensional
Bessel process. For any $\d_2\in (0, \, {1\over 2})$ and all
large~$t$,
$$
\P \left\{\left|\frac{1}{\log t} \int_0^{t}
\frac{\textnormal{d}s}{R^2(s)} - \frac{1}{d-2}\right|> {1\over (\log
t)^{(1/2)-\d_2}} \right\} \le \exp \left( - c_{20}\, (\log
t)^{2\d_2} \right).
$$

\end{lemma}

\begin{lemma}\label{lemmaApproximationU}
If $\d_1>0$ is small enough, then for all large $v$, $\P(\EG^c)\leq
v^{-1/4+5\d_1}$, where
\begin{equation}
\EG:=\{\tau_{\beta}[(1-v^{-\d_1})\l v ]\leq U(v)\leq
\tau_{\beta}[(1+v^{-\d_1})\l v]\}.\label{eqLemmaApproxUbis}
\end{equation}

\end{lemma}

\begin{lemma}
 \label{remarquenegatif}
 Let $\k>0$ and define
$
     H_-(+\infty)
   :=\lim_{r\to+\infty}H_-(r).
$
 There exists
 $c_{21}>0$ such that for all large
 $z$,
 \begin{eqnarray}
     \P(H_-(+\infty)>z)
  & \leq & c_{21}[(\log z)/z]^{\k/(\k+2)}.
    \label{eqHmoins}
 \end{eqnarray}
\end{lemma}

\begin{lemma}\label{lemmaApproxJCasKinf1}
Let $(\b(t),\ t\geq 0)$ be a Brownian motion, and
let $\l=4(1+\k)$. We define
\begin{equation}
J_{\b}(\k, t,\l)  :=  \int_0^1
y(1-y)^{\k-2}L_\b(\tau_\b(\l), {S(y)\over t} )\text{d}y, \qquad
0<\k\leq 1,\ t\geq 0. \label{e1p18}
\end{equation}
Let $0<d<1$ and let $0<\e<1$. There exists $c_{22}>0$ such that for
$t$ large enough, on an event $\ER$ of probability greater than
$1-c_{22}/t^d$,

 \noindent {\bf (i)} Case $0<\k<1$: (recall
$K_\b({\k})$ from \eqref{e3p18})
\begin{equation}
\label{eqLemmaJcasKinf1} (1-\e)K_\b(\k)- c_{23} t^{1-1/\k}\leq
{\k}^{2-1/\k}t^{1-1/\k}J_\b(\k,t,\l)\leq (1+\e)K_\b(\k)+ c_{23}
t^{1-1/\k}.
\end{equation}

\noindent {\bf (ii)} Case $\k=1$: (recall $C_\b$ from \eqref{e2p18})
\begin{equation}
\label{eqLemmaJcasKeg1} (1-\e)[C_\b + 8\log t]\leq J_\b(1,t,
8)\leq (1+\e)[C_\b + 8\log t].
\end{equation}
\end{lemma}

\medskip

By admitting Lemmas
\ref{lemmaApproxTy}--\ref{lemmaApproxJCasKinf1}, we can now
complete the proof of Lemma \ref{lemmaApproxLxLi}.

\medskip

\noindent {\bf Proof of Lemma \ref{lemmaApproxLxLi}.}
Notice that
\begin{equation}
    S(y) {\sim}_{y\to 1}
    (1-y)^{-\k}/\k, \qquad
y/(1-y)  {\sim}_{y\to 1} [\k S(y)]^{1/\k}. \label{equivalent}
\end{equation}

We first look for an estimate of $U[\g(r)-T_Y(\a_\k)]$. Since
$A_{\infty}\el 2/\g_\k$, where $\g_\k$ is a gamma variable of
parameter $\k$, and  $A(F(r))\le A_{\infty}$, we have $
\P[A(F(r))>r^{2/\k}] \le \frac{2^\k r^{-2}}{\k\Gamma(\k)}.$ On the
other hand, by definition, $A(F(r)) = A_{\infty}-\delta(r) =
A_{\infty}- e^{-\k r/2}$ (see \eqref{e1p5}). Hence,
$$
\P \left[ A(F(r)) <{1/( 2\log r)} \right] \le \P \left[ {2/ \g_\k} <
{1/( 2\log r)} + \delta(r) \right] \le r^{-2}/\Gamma(\k) .
$$


\noindent Recall that $\g(r)=\ty[\d(r)^{-1}A(F(r))]$, see
\eqref{e10p7}. Thus, for large $r$,
$$
\P \left\{ \ty[\exp(\k r/2 - 2\log\log r)] \le \g(r)\leq \ty[
\exp(\k r/2 +(2/\k) \log r)] \right\} \ge 1- {c_{24} r^{-2}}.
$$

\noindent By definition, $\ty (u) = \int_0^u
\frac{\text{d}s}{\widetilde{R}_2^2(s) + \rb^2(s+1)}$. Since
$(\widetilde{R}_2^2(t)+\rb^2(t+1),\ t\ge 0)$ is a
$(4+2\k)$--dimensional squared Bessel process starting from
$\Rt_{2+2\k}^2(1)$, it follows from Lemma \ref{lemmaApproxTy} that
there exist constants $\delta_3\in (0, \frac{1}{2})$, $c_{25}>0$ and
$c_{26}>0$, such that for large $r$, with $\l = 4(1+\k)$ as before,
\begin{equation}
    \label{gamma}
    \P \left\{\k r/\l-c_{25}r^{1/2+\d_2}\leq \g(r)\leq
    \k r/\l+c_{25}r^{1/2+\d_2}\right\} \ge 1- {c_{26}
     r^{-2}}.
\end{equation}

To study the behaviour of $T_Y(\a_\k)$, we notice that $Y$ satisfies
\eqref{eqDefJacobi} with $d_1=2$ and $d_2=2+2\k$. By the
Dubins--Schwarz theorem, there exists a Brownian motion $(\Bhat(t),\
t\geq 0)$ such that $ Y(t)=\Bhat\left(4\int_0^t Y(s)(1-Y(s))\text{d}
s\right)+\int_0^t [2-(4+2\k)Y(s)]\text{d} s$ for $t\geq 0. $ Recall
that $\a_\k=1/(4+2\k)$. Let $t\ge 2\a_\k$. We have, on the event $\{
T_Y(\a_\k)\ge t \}$,
$$
\inf_{0\le s\le 4 t}\Bhat(s)\leq \Bhat\big(4\int_0^{t}
Y(s)(1-Y(s))\text{d} s\big)\leq \a_\k-t\leq -\frac{t}{2},
$$

\noindent since $Y(s)\in(0,1)$ for any $s\ge 0$. As a consequence,
for $t\ge 2\a_\k$,
\begin{equation}
    \label{eqLemmaHittingTimeYConclusion}
    \P(T_Y(\a_\k)>t) \le \P[\inf\nolimits_{0\le s\le 4t}
    \Bhat(s)\le -t/2]\le 2 \exp [-t/32].
\end{equation}

\noindent In particular, $\P[T_Y(\a_\k)> 64\log r ] \le {2\over
r^2}$ for large $r$. Plugging this into (\ref{gamma}), and
introducing ${\underline \g} = {\underline \g} (r) :=\frac{\k}{\l}
r- 2c_{25} r^{1/2+\d_2}$ and ${\overline \g} = {\overline \g} (r)
:=\frac{\k}{\l} r+ c_{25} r^{1/2+\d_2}$ yields that for large $r$,
$$
\P \left\{ U({\underline \g})\leq U[\g(r)-T_Y(\a_\k)] \leq
U({\overline \g}) \right\} \ge 1- {c_{27}  r^{-2}}.
$$

\noindent By Lemma \ref{lemmaApproximationU}, for small $\d_1>0$
and all large $r$,
$$
\P\left\{ \tau_\beta \left[ (1- {{\underline \g}^{-\d_1}} ) \l
{\underline \g}\, \right] \le U[\g(r)-T_Y(\a_\k)] \le \tau_\beta
\left[ (1+ { {\overline \g}^{-\d_1}} ) \l {\overline \g} \right]
\right\} \ge 1- { r^{-c_{28}}}.
$$

\noindent We choose $\d_1$ so small that $\d_1<1/2-\d_2$. Then for
large $r$, we have $(1- {1\over {\underline \g}^{\d_1}} ) \l
{\underline \g} \ge [1-2({\l \over \k})^{\d_1}r^{-\d_1}] \k r=\l
t_-(r) $, and $(1 + {1\over {\overline \g}^{\d_1}} ) \l {\overline
\g} \le [1+2({\l \over \k})^{\d_1}r^{-\d_1}] \k r=\l t_+(r)$\ ¨(see
\eqref{e1p22}). Hence,
\begin{equation}
    \P\left\{ \tau_\beta[\l t_-(r)] \le U[\g(r)-T_Y(\a_\k)]
    \le \tau_\beta[\l t_+(r)] \right\} \ge 1- {
    r^{-c_{28}}}.
    \label{eqEncadrementUti}
\end{equation}


As in Hu et al. (1999), p.~3923), \eqref{e1p7} leads to $ H_0(r) =
4\int_0^1 x (1-x)^{\k-2} L_{\beta} [U(\g(r)- T_Y(\a_\k)), S(x)]
\text{d} x$.
 By \eqref{eqEncadrementUti},  $
    \P \left[ I_-'(r) \le H_0 (r) \le I_+'(r) \right] \ge
    1-  r^{-c_{28}} ,
$ for large $r$, where
\begin{equation*}
    I_\pm'(r)
  :=  4\int_0^1  x (1-x)^{\k-2} L_{\beta}\{ \tau_{\beta} [\l
    t_\pm(r)], S(x) \} \text{d}x
     =  4 t_\pm(r) J_{\beta_{t_\pm(r)}} [\k,t_\pm(r), \l],
\end{equation*}

\noindent and, as before, $t_\pm(r)= [1\pm 2 ({\l \over
\k})^{\delta_1}r^{-\d_1} ] {\k \over \l} r$, $\b_v(s)=\b(v^2 s)/v$
and $J_{\b}$ is defined in \eqref{e1p18}.

Applying Lemma \ref{lemmaApproxJCasKinf1} to $d=1/2$
yields that, for large $r$, recalling $\Ib_\pm (r)$ from
(\ref{eqI''CasKinf1}),
\begin{equation}
    \P \{ (1-\e)\Ib_- (r)\leq H_0(r)\leq (1+\e)\Ib_+ (r) \} \ge
    1- {r^{-c_{29}}}.
    \label{eqI''CasGeneral}
\end{equation}

In the case $0<\k<1$, \eqref{eqLemmaHittingTimeYConclusion} and
\eqref{eqMajorAvantAlphaK} give $\P[\hbar(r)\le c_{30}\log r]\geq
1-2r^{-2}$ for some $c_{30}$ and all large $r$.
On the other hand, by Lemma \ref{remarquenegatif}, $\P[H_-(F(r))\le
\e r] \ge \P[H_-(+\infty)\le \e r] \ge 1- {c_{31}\over
r^{(1-\d_1)\k/(\k+2)}}$, for all large $r$. Consequently, by
\eqref{eqI''CasGeneral} and \eqref{eqMajorAvantAlphaK}, for large
$r$,
$$
\P[ (1-\e)\Ib_-(r)\leq H(F(r))\leq (1+\e)\Ib_+(r)+ (4\e\l/\k)
t_+(r)] \geq 1-{r^{-c_{32}}}.
$$

\noindent This proves Lemma \ref{lemmaApproxLxLi} in the case $0<\k<1$.

Now we turn to the case $\k=1$. We again have
$\P[H_-(F(r))+\hbar(r)\leq 2\e r]\ge 1- r^{-c_{33}}$ (for large
$r$). By Biane and Yor
(1987)
$C_{\b_{t_\pm(r)}} =^{\mathcal{L}} \pi C_8^{ca}/2+c_{34}$,
where $c_{34}>0$. Hence, by \eqref{ProbaSamoTaqqu},
$\P[C_{\b_{t_\pm(r)}}>-\pi\log r]\geq 1-r^{-2}$. Thus
\eqref{eqI''CasKinf1} gives
$ \P\{\Ib_+(r)\geq 16t_+(r)\log r\}\geq 1-r^{-2}$. This, together
with \eqref{eqI''CasGeneral}, yields that, for large $r$,
$$
\P[ (1-\e)\Ib_-(r)\le H(F(r))\le (1+2\e)\Ib_+(r)]\ge 1-{
r^{-c_{35}}} .
$$

\noindent This proves Lemma \ref{lemmaApproxLxLi} in the
case $\k=1$. \hfill$\Box$

\medskip

The rest of the section is devoted to the proof of Lemmas
\ref{lemmaApproxTy}--\ref{lemmaApproxJCasKinf1}.


\medskip

\noindent {\bf Proof of Lemma \ref{lemmaApproxTy}.} Let $d>4$ and $R_0\el\widetilde{R}_{d-2}(1)$, where
$\widetilde{R}$ is a $(d-2)$--dimensional Bessel process. We
consider a $d$--dimensional Bessel process $R$, starting from
$R_0$. Let $\theta(t):=\int_0^t R^{-2}(s) \text{d}s$. It\^o's
formula gives $\log R(t)=\log R_0
+M(t)+\frac{d-2}{2}\theta(t)$, where $M(t):=\int_0^t
R(s)^{-1}\text{d}\bhat(s)$ and $(\bhat(t),\ t\geq 0)$ is a
Brownian motion. By the Dubins--Schwarz theorem, there exists
a Brownian motion $(\Bt(t),\ t\geq 0)$ such that
\begin{equation}
    \label{ito}
    (d-2) \theta(t)/2 = \log R(t) - \log R_0
    - \Bt(\theta(t)), \qquad t\geq 0.
\end{equation}

Let $\delta_3 \in (0, {1\over 2})$, and let $x=x(t) := {d-2\over
6} {1\over (\log t)^{(1/2)-\delta_3}}$.
For large $t$, we have
\begin{equation}
    \label{eqTpuissanceEpsilon1}
    \P\left(\left |\frac{\log R_0}{\log t}\right|> x\right)
    \leq \P\left( \frac{\log R_0}{\log t}> x \right)
    +\P\left(
\frac{\log R_0}{\log t} < -x \right)
    \le \exp\left( - (1-\varepsilon) {t^{2x}\over 2}
    \right) + c_{36} \, t^{-\frac{x}{d-2}} .
\end{equation}

Let $n:= \lceil d\rceil$ be the smallest integer such that $n\ge d$.
Since an $n$-dimensional Bessel process can be realized as the
Euclidean modulus of an $\R^n$-valued Brownian motion, it follows
from the triangular inequality that $ R(t)\leq_{\mathcal{L}} R_0
+\widehat{R}_n(t)$,
 where $(\widehat{R}_n(t),\ t\geq 0)$ is an
$n$-dimensional Bessel process starting from $0$. Consequently,
$ \P\left(R(t)>t^{(1/2)+x}\right) \le \exp\left( - (1- \e) {t^{2x}/
4}\right)$ for large $t$, and $ \P\left(R(t)<t^{(1/2)-x} \right)
\leq c_{36} \, t^{-x/d}$. Therefore, for large $t$,
\begin{equation}
    \label{eqTpuissanceEpsilon11}
    \P\left(\left |(\log R(t))/\log t-1/2\right|> x \right)
    \le \exp\left( - (1- \e) {t^{2x}/ 4}\right) +
    c_{36} \, t^{-x/d} .
\end{equation}

Define $    \EY:=\{ \sup\nolimits_{0\le s \leq 2(\log t)/(d-2)}
|\Bt(s)| \leq
    x\log t\}
$ and
\begin{eqnarray*}
    \EK:= \left\{ \left|\frac{\log R(t)}{\log t} - \frac{1}{2}
    \right| \le x \right\} \cap \left\{ \left|
    \frac{\log R_0}{\log t} \right| \le x \right\},
    &\quad \ET :=& \left\{ \frac{d-2}{2} \theta(t) < \log t\right\}.
\end{eqnarray*}

By \eqref{eqTpuissanceEpsilon1} and
\eqref{eqTpuissanceEpsilon11}, we have, for large $t$,
\begin{equation}
\P\left(\EK^c\right) \le 2\exp\left( - (1- \e) {t^{2x}/ 4}\right) +
c_{37} \, t^{-x/d} .\label{Mand1}
\end{equation}

We now estimate $\P(\EK\cap\ET^c)$. We first observe that  $ \left|
\Bt(\theta(t)) + (d-2) \theta(t)/2- (\log t)/2 \right|\le 2x\log t$
on $\EK$, by \eqref{ito}. We claim that $\EK\cap\ET^c \subset
\{|\Bt(\theta(t))|>\frac{d-2}{6} \theta(t)\}$ for large $t$. Indeed,
on the event $\EK\cap\ET^c \cap \{|\Bt(\theta(t))| \le \frac{d-2}{6}
\theta(t)\}$,
$$
(d-2)\theta(t)/2 \le \left(2x+1/2 \right) \log t- \Bt(\theta(t)) \le
\left(2x+1/2\right) \log t + (d-2) \theta(t)/6 ,
$$

\noindent which implies $\frac{d-2}{2}\theta(t) \le
(\frac{3}{4}+3x)\log t$. This, for large $t$, contradicts
$\frac{d-2}{2}\theta(t)>\log t$ on $\ET^c$. Therefore, $\EK\cap\ET^c
\subset \{|\Bt(\theta(t))|>\frac{d-2}{6} \theta(t)\}$ holds for all
large $t$, from which it follows that
\begin{equation*}
    \P(\EK\cap\ET^c)
 \le \P\left( \, \sup\nolimits_{s\ge 2(\log t)/(d-2)}
    |\Bt(s)|/s > (d-2)/6 \right)
 \le \exp\left[ - (1-\e) {(d-2)(\log t) /36}\right].
\end{equation*}

\noindent Since $\P(\EY^c)\le \exp[- (1-\e) {d-2\over 4} x^2 \log
t]$ (for large $t$), this and \eqref{Mand1}  give  for large $t$,
\begin{eqnarray*}
    \P(\EK^c \cup \ET^c \cup \EY^c)
 \le  \P(\EK^c)+\P(\EK\cap\ET^c)+P(\EK\cap\ET\cap\EY^c)
    \le \exp(- c_{38} \, x^2 \log t).
\end{eqnarray*}

\noindent Since $\EK\cap \ET \cap \EY \subset \{ | {\theta (t)\over
\log t} - {1\over d-2}| \le {6x\over d-2}\}$, this completes the
proof of Lemma \ref{lemmaApproxTy}.
\hfill$\Box$


\medskip

\noindent {\bf Proof of Lemma \ref{lemmaApproximationU}.} Let $v>0$. Recall that for $x>0$, $\b_v(x) = (1/v)\b(v^2x)$, and
$v^2 \tau_{\beta_v}(x)=\tau_{\beta}(x v)$ a.s. Then,
\begin{equation}
\EG=\left\{\tau_{\beta_v}[(1-v^{-\d_1})\l]\leq U(v)/v^2\leq
\tau_{\beta_v}[(1+v^{-\d_1})\l]\right\}. \label{eqLemmaApproxU}
\end{equation}

\noindent For $\d_1>0$, define $ \EGb  := \{\sup\nolimits_{0\leq
s\leq \tau_{\beta_v}(2\l)}|\e_1(v,s)|<v^{-\d_1}\}$, where
$\e_1 =  \e_1(v,s)  :=  \frac{1}{4}\int_0^1
(1-x)^{\k}\left[L_{\beta_v}(s, {S(x)/
v})-L_{\beta_v}(s,0)\right]\text{d}x$ for $s\geq 0$.
  By Hu et al. (1999)
  p.~3924), $\EGb\subset \EG$. Thus it remains to prove that for
$\d_1$ small enough, $\P(\EGb^c)\leq 1/v^{1/4-5\d_1}$ for large $v$.
Notice that for $s\geq 0$,
\begin{eqnarray}
|\e_1| & \leq & \left(
\int_{\{S(x)>\sqrt{v}\}}+\int_{\{S(x)<-\sqrt{v}\}}+
\int_{\{|S(x)|\leq\sqrt{v}\}} \right)
\frac{(1-x)^{\k}}{4}\left|L_{\beta_v}(s, {S(x)/
v})-L_{\beta_v}(s,0)\right|\text{d}x
\nonumber\\
&=:& \e_2(v,s)+\e_3(v,s)+\e_4(v,s)\label{voila2}.
\end{eqnarray}

By \eqref{equivalent}, we have, for all large $v$ (and all $s\ge 0$)
\begin{eqnarray*}
    \sup_{0\leq s \leq \tau_{\beta_v}(2\l)}\e_2(v,s)
    & \leq &
    \frac{1}{4}\int_{1-\left(\frac{2}{\k\sqrt{v}}\right)^{1/\k}}^1 \
    (1-x)^{\k} \sup_{0\leq s \leq \tau_{\beta_v}(2\l)}
    \sup\limits_{u\geq
0}\left[L_{\beta_v}(s,u)+L_{\beta_v}(s,0)\right]\text{d}x\\
    \sup\nolimits_{u\geq
0}\left[L_{\beta_v}(\tau_{\beta_v}(2\l),u)+2\l\right].
\end{eqnarray*}

\noindent By the second Ray--Knight theorem,
$Z:=(L_{\beta_v}(\tau_{\beta_v}(2\l),u), \ u\geq 0)$ is a
$0$--dimensional squared Bessel process starting from $2\l$.
Hence, for large $v$,
\begin{equation}
\P[\sup\nolimits_{0\leq s \leq \tau_{\beta_v}(2\l)}\e_2(v,s)\geq
\left[2/(\k\sqrt{v})\right]^{\frac{1}{\k}+1}(\sqrt{v}+2\l)] \leq
\P\left(\sup\nolimits_{u\geq 0}Z(u)\geq \sqrt{v}\right)
 =2\l/\sqrt{v}.\label{eqInegEpsilon2}
\end{equation}

Similarly (this time, using $S(x)\sim\log x$, $x\to 0$), we
have, for large $v$,
\begin{equation}
\label{eqInegEpsilon3}
    \P[\sup\nolimits_{0\leq s \leq
    \tau_{\beta_v}(2\l)}\e_3(v,s)\geq \exp\left(-\sqrt{v}/2\right)
    (\sqrt{v}+2\l)] \leq 2\l/\sqrt{v}.
\end{equation}

To estimate $\e_4(v,s)$, we note that
$
\e_4(v,s)
 \leq  \sup\nolimits_{|u| \le 1/\sqrt{v}}
\left|L_{\beta_v}(s,u)-L_{\beta_v}(s,0)\right|.
\label{eqInegEpsilon4}
$

\noindent Let $\e\in(0,1/2)$, $t_v>0$, $\g\geq 1$ and define
$(\beta_v)_{t_v}^*:=\sup_{0\leq s\leq t_v}|\beta_v(s)|$. Applying
Barlow and Yor (1982)
(ii)) to the
continuous martingale $\beta_v(.\wedge t_v)$, we see that for
some constant $C_{\g,\e}>0$,
$$
\|\sup\nolimits_{0\leq s\leq t_v,a\neq b}
|L_{\beta_v}(s,b)-L_{\beta_v}(s,a)|/|b-a|^{1/2-\e}\|_{\gamma} \leq
C_{\g,\e}\|[(\beta_v)_{t_v}^*]^{1/2+\e}\|_{\g}.
$$
Then, by Chebyshev's inequality, for $\a>0$,
\begin{equation}
    \P\big(\sup\nolimits_{0\leq s\leq t_v, \, a\neq b}
    \frac{|L_{\beta_v}(s,b)
    -L_{\beta_v}(s,a)|}{|b-a|^{1/2-\e}}\geq
    \a \big)
\le \frac{(\sqrt{t_v})^{(1/2+\e)\g}}{\a^{\g}}
    \left[C_{\g,\e}\| [(\beta_v)_{1}^*]^{1/2+\e}
    \|_{\g}\right]^\g.
    \label{eqApplicationBarlowYor}
\end{equation}

On $
\EH  :=  \{ \sup\nolimits_{0\leq s\leq \tau_{\beta_v}(2\l),a\neq b}
|L_{\beta_v}(s,b)-L_{\beta_v}(s,a)|/|b-a|^{1/2-\e} \leq
v^{\frac{1}{2}\left(\frac{1}{2}-2\e\right)} \},
$ we have
\begin{eqnarray}
\sup\nolimits_{0\leq s \leq \tau_{\beta_v}(2\l)}\e_4(v,s) & \leq &
v^{\frac{1}{2}(-\frac{1}{2}+\e)}v^{\frac{1}{2}
\left(\frac{1}{2}-2\e\right)}=v^{-\e/2}. \label{eqInegEpsilon4surEH}
\end{eqnarray}
We choose $\g:=2$ and $t_v:=v^{\frac{1/4-\e}{1/2+\e}}$ to see that
for all large $v$ (if $\e$ is small enough),
\begin{eqnarray*}
    \P(\EH(v)^c)
    \leq  \P\left(\tau_{\beta_v}(2\l)>t_v\right)+
        (\sqrt{t_v})^{(1/2+\e)\g}
    \left[C_{\g,\e}\| [(\beta_v)_{1}^*]^{1/2+\e}
    \|_{\g}\right]^\g
    (v^{1/4-\e})^{-\g}
    \leq v^{-1/4+2\e}/2.
\end{eqnarray*}

\noindent Combining this with  \eqref{voila2},
\eqref{eqInegEpsilon2}, \eqref{eqInegEpsilon3} and
\eqref{eqInegEpsilon4surEH}, we obtain that, for $\e>0$ small
enough,
$$
\P(\sup\nolimits_{0\leq s\leq \tau_{\beta_v}(2\l)}|\e_1(v,s)|\geq
2v^{-\e/2} )\leq v^{-1/4+2\e}.
$$

\noindent This gives, with the choice of $\d_1:=2\e/5$,
$\P(\EGb^c)\leq v^{-1/4+5\d_1}$ (for large $v$).\hfill$\Box$


\medskip

\noindent {\bf Proof of Lemma \ref{remarquenegatif}.} For $a>0$, $\a>0$ and $b>0$, let
\begin{eqnarray*}
\EM  := \{\sup\nolimits_{x<0} e^{-\wk(x)}\leq a \}, \quad \EN :=
\left\{A_{\infty}\leq\a\right\},\quad
 \EO  :=\left\{\sup\nolimits_{y<0} \lB[\s_B(\a),y]\leq
 b\right\},&&\\[1mm]
\lneg  :=  \sup\limits_{r\geq 0} \ \sup\limits_{x<0}
\left\{e^{-\wk(x)}\lB[\s_B(A(r)),A(x)]\right\}
 \leq  \left(\sup\limits_{x<0} e^{-\wk(x)}\right)
\sup\limits_{y<0} \lB[\s_B(A_{\infty}),y].&&
\end{eqnarray*}

It follows from the second Ray--Knight theorem that $\P(\EO^c)\leq
c_{39} \a/b$. Now, let $a:=z^{\frac{1}{\k+2}}$,
$\a:=z^{\frac{1}{\k+2}}$ and $b:=z^{\frac{\k+1}{\k+2}}$. Notice that
$\lneg\leq z$ on $\EM\cap\EN\cap\EO$, and recall $A_{\infty}\el
2/\g_{\k}$,
where $\g_{\k}$ is a gamma variable of parameter $\k$. We have for
$z$ large enough,
\begin{eqnarray}
\P(\lneg>z)  \leq \P(\EM^c)+\P(\EN^c)+\P(\EO^c)
 \leq a^{-\k}
+c_{39}\left(2/\a\right)^\k +c_{39}\a/b  \leq c_{40}
z^{-\frac{\k}{\k+2}}.\ \ \label{eqprobasEMNO}
\end{eqnarray}

Define for $c>0$, $\EJ  :=  \{\min\nolimits_{0\leq s\leq
\s_B(A_{\infty})}B(s)> -A_{\infty}z^{\frac{\k+1}{\k+2}}\}$, $ \ES
:=  \{|A^{-1}(-z)|\leq c\log z\}. $
 On $\EM\cap\dots\cap\ES$, $ H_-(+\infty)
 \leq  \lim_{r\to+\infty} \int_{A^{-1}\left(\min_{0\leq s\leq
    \s_B(A(r))}B(s)\right)}^0
    \lneg\text{d}x
$
 for $r\geq 0$. Hence,
\begin{eqnarray}
H_-(+\infty)  \leq  \left|A^{-1}\left(\min\nolimits_{0\leq s\leq
    \s_B(A_{\infty})}B(s)\right)\right|\lneg
 \leq  |A^{-1}(-z)|\lneg
  \leq  c z\log z.\label{exp47}
\end{eqnarray}

\noindent Moreover, for $c>2/\k$, $\e>0$, and $z$ large enough,
\begin{eqnarray}
    \P(\ES^c)
    = \P\left(z>\int_0^{c\log z}e^{W(u)+\k u/2}\text{d}u\right)
    & \leq & \P\left[z>\exp\left(\inf_{0\leq u\leq c\log z}W(u)\right)\frac{2}{\k}(z^{\k
    c/2}-1)\right]
\nonumber\\
       & \leq & 2z^{-\frac{1}{2c}\left(\frac{\k c}{2}-1-\e\right)^2}.
\label{AsterixCielsurlatete}
\end{eqnarray}

Since $B$ is independent of $A_{\infty}$, we have $ \P\left(\EJ^c
|A_{\infty}\right)=A_{\infty}/[A_{\infty}+A_{\infty}z^{\frac{\k+1}{\k+2}}
]\leq z^{-\frac{\k+1}{\k+2}}. $ Choosing $c$ large enough, this,
together with \eqref{eqprobasEMNO}, \eqref{exp47} and
\eqref{AsterixCielsurlatete}, gives \eqref{eqHmoins}. \hfill$\Box$


\medskip

\noindent {\bf Proof of Lemma \ref{lemmaApproxJCasKinf1}.} Assume
$0<\k\leq 1$. Consider a Brownian motion $\b$, a small constant
$\e>0$, and $0<d<1$. Recall
$S(y)=\int_{\a_{\k}}^y\frac{\text{d}x}{x(1-x)^{1+\k}} $ and notice
that $ 1-S^{-1}(u){\sim}_{u\to+\infty}(\k u)^{-1/\k}. $ Therefore,
there exists $x_\e>0$ such that for all $u\geq x_\e$ ,
$[1-S^{-1}(u)]^{2\k-1}/(\k u)^{1/\k-2}\in (1-\e,1+\e)$  and
$S^{-1}(u)\geq (1-\e)$.

 Let $g(t):=t^{\e-1}$, and write
\begin{eqnarray*}
 & &    J_{\b}(\k,t,\l)
\\
    & = & \left(\int_{\{S(y)\leq -t g(t)\}}+\int_{\{-t g(t)< S(y)\leq
    0\}} +\int_{\{0< S(y)\leq x_\e\}}+\int_{\{x_\e<S(y)\}}\right)
    \frac{y}{(1-y)^{2-\k}}L_\b(\tau_\b(\l),{S(y)\over t})\text{d}y
\\
    & := & J_1+J_2+J_3+J_4.
\end{eqnarray*}

Since $S(x)\sim\log x$, $x\to 0$, we have for large $t$,  $ J_1 \leq
\exp\left(-\frac{t g(t) }{2}\right)\left(\sup_{s\geq 0}
Z(s)\right)$, where $Z$ is a $0$--dimensional squared Bessel process
starting from $\l$ (by the second Ray--Knight theorem). Hence, we
get $ \P\left[J_1 \geq e^{-t^\e/2}t^d \right] \leq \l/t^d. $

\noindent Fix a large constant $\g>0$, and define
\begin{eqnarray*}
\EP  :=  \{\tau_\b(\l)\leq t^{2d}\},\qquad \EQ  :=
\left\{\sup\nolimits_{0\leq s\leq t^{2d}, \, a\neq b
}{|L_\b(s,b)-L_\b(s,a)|}/{|b-a|^{1/2-\e}}\leq
t^{d(1/2+\e+1/\g)}\right\}.
\end{eqnarray*}

Recall that $S(\a_\k)=0$. On the event $\EP\cap\EQ$ and for
all large $t$,
\begin{eqnarray*}
    \k^{2-1/\k}t^{1-1/\k}J_3
    & \leq & \k^{2-1/\k}
    t^{1-1/\k}\sup\limits_{0\leq x\leq x_\e/t} L_\b(\tau_\b(\l),x)
    \int_{\a_k}^{S^{-1}(x_\e)}y(1-y)^{\k-2}\text{d}y
\\
& \leq & c_{41}
t^{1-1/\k}[\l+t^{d(1/2+\e+1/\g)}(x_\e/t)^{\frac{1}{2}-\e}]
  \leq  2\l c_{41}t^{1-1/\k} ,
\end{eqnarray*}
\begin{equation*}
    J_2
    \leq \sup\nolimits_{-g(t)\leq s\leq 0} L_\b(\tau_\b(\l),s)
    \int_0^{\a_k}y(1-y)^{\k-2}\text{d}y
\leq  c_{42}\left[\l+ t^{d(1/2+\e+1/\g)}
    (t^{\e-1})^{\frac{1}{2}-\e}\right]  \leq 2 c_{42}.
    \end{equation*}

\noindent Since $\P(\EP^c)\leq c_{43}/t^d$ and $\P(\EQ^c)\leq
c_{44}/t^d$ (see \eqref{eqApplicationBarlowYor}), we obtain, for
large $t$, $
    \label{eqApproxJ3b}
    \P\left( J_3 \le c_{45} \right) \ge 1 - {c_{46}
    /t^d}
$ and
$
    \P\left( J_2 \le 2 c_{42} \right) \ge 1 - {c_{47}
    /t^d}.
$

Now, we write $ J_4
 =  \k^{1/\k-2}\ t^{1/\k-1}\int_{x_\e/t}^{+\infty}\left(S^{-1}(t
x)\right)^2\frac{\left(1-S^{-1}(t x)\right)^{2\k-1}}{(\k
t)^{1/\k-2}}L_\b(\tau_\b(\l),x)\text{d}x. $  Therefore
\begin{eqnarray}
(1-\e)^3\int_{x_\e/t}^{+\infty}x^{\frac{1}{\k}-2}L_\b(\tau_\b(\l),x)\text{d}x
\leq   \k^{2-\frac{1}{\k}}t^{1-\frac{1}{\k}}J_4
  \leq
(1+\e)\int_{x_\e/t}^{+\infty}x^{\frac{1}{\k}-2}L_\b(\tau_\b(\l),x)\text{d}x.
\ \ \label{eqApproxJ3a}
\end{eqnarray}

We first assume $0<\k<1$. On $\EP\cap\EQ$, for large $t$, we have
$\int_0^{x_\e/t}x^{1/\k-2}L_\b(\tau_\b(\l),x)\text{d}x\leq c_{48}
t^{1-1/\k}$. Recall $K_{\b}(\k)$ from \eqref{e3p18}. By
\eqref{eqApproxJ3a}, for large~$t$,
$$
\P\left[ (1-\e)^3K_\b(\k)-(1-\e)^3c_{48} t^{1-1/\k} \le
\k^{2-1/\k}t^{1-1/\k} J_4 \le (1+\e)K_\b(\k) \right] \ge 1- {c_{49}/
t^d},
$$

Since $J_{\b}(\k,t,\l) = J_1+J_2+J_3+J_4$, we get
$$
\P\left\{ (1-\e)^3K_\b(\k)-c_{23}t^{1-1/\k}\leq
{\k}^{2-1/\k}t^{1-1/\k}J_\b(\k,t,\l)\leq
(1+\e)K_\b(\k)+c_{23}t^{1-1/\k} \right\}\geq 1- \frac{c_{50}}{t^d},
$$

\noindent proving the lemma in the case $0<\k <1$.

We now assume $\k=1$. By the definition of $C_\b$ (see
\eqref{e2p18}),
\begin{eqnarray*}
\int_{x_\e/t}^{\infty}\frac{L_\b(\tau_{\beta}(8),x)}{x}\text{d}x
 & = & C_\b -\int_0^{x_\e/t}
\frac{L_\b(\tau_{\beta}(8),x)-8}{x}\text{d}x+8\log t-8\log
    x_\e.
\end{eqnarray*}

\noindent On $\EP\cap\EQ$, for large $t$, $ \int_0^{x_e/t}
\frac{|L_\b(\tau_\b(8),x)-8|}{x}\text{d}x  \leq
 \int_0^{x_e/t} \frac{t^{d(1/2+\e+1/\g)}x^{1/2-\e}}{x}\text{d}x
\leq\e$. As in \eqref{ProbaSamoTaqqu}, $\P(C_\b+8\log t <\log t
)\leq r^{-7}$. Therefore, by (\ref{eqApproxJ3a}), we have, for large
$t$,
$$
\P\left\{ (1-\e)^4 [C_\b +8\log t]\le J_4 \le (1+\e)^2 [C_\b +8\log
t] \right\} \ge 1- {c_{51}/ t^d} .
$$

Since $J_{\b}(1,t,8) = J_1+J_2+J_3+J_4$, this yields that for
large~$t$,
$$
\hspace{1,9cm} \P\left\{ (1-\e)^4 [C_\b +8\log t] \le J_\b(1,t,8)\le
(1+\e)^3 [C_\b +8\log t] \right\} \ge 1- {c_{52}/
t^d}.\hspace{1,9cm}\Box
$$

\noindent {\bf  Acknowledgements.} I would like to thank Zhan Shi
for many helpful discussions.


\addcontentsline{toc}{section}{Bibliography}

\bigskip
\noindent{\bf References :}

{\baselineskip=14pt

\noindent
    Barlow, M. T. and Yor, M. (1982),    Semimartingale inequalities via the
    {G}arsia-{R}odemich-{R}umsey lemma, and applications to
    local times.
    {\it J. Funct. Anal.} {\bf 49},  198--229.

\medskip\noindent
    Bertoin, J. (1996),
    {\it L\'evy Processes}.
    Cambridge University Press, Cambridge.

\medskip\noindent
    Biane, Ph. and Yor, M. (1987),
    Valeurs principales associ\'ees aux temps locaux browniens.
    {\it Bull. Sci. Math.} {\bf 111},  23--101.

\medskip\noindent
    Borodin, A. N. and Salminen, P. (2002),
    {\it Handbook of {B}rownian Motion---Facts and Formulae}.
    {Birkh\"auser}, Boston.

\medskip\noindent
    Brox, Th. (1986),
    A one-dimensional diffusion process in a {W}iener medium.
    {\it Ann. Probab.} {\bf 14},  1206--1218.

\medskip\noindent
    Devulder A. (2005),
    The maximum local time of a diffusion
    in a brownian potential with drift.
    Preprint.

\medskip\noindent
    Dufresne, D. (2000),
    Laguerre series for {A}sian and other options.
    {\it Math. Finance} {\bf 10}, 407--428.

\medskip\noindent
 Hu, Y., Shi, Z. and Yor, M. (1999),
    Rates of convergence of diffusions with drifted {B}rownian
    potentials.
    {\it Trans. Amer. Math. Soc.} {\bf 351},  3915--3934.

\medskip\noindent
    Kawazu, K. and Tanaka, H. (1997),
    A diffusion process in a Brownian environment with drift.
    {\it J. Math. Soc. Japan} {\bf 49},  189--211.
\medskip\noindent
    Samorodnitsky, G. and Taqqu, M. S. (1994),
    {\it Stable Non-{G}aussian Random Processes}.
    Chapman \& Hall,  New York.

\medskip\noindent
    Schumacher, S. (1985),
    Diffusions with random coefficients.
    {\it Contemp. Math.} {\bf 41},  351--356.

\medskip\noindent
    Warren, J. and Yor, M. (1997),
    Skew products involving Bessel and Jacobi processes. Technical
    report, Statistics group, University of Bath.

\medskip\noindent
    Zeitouni, O. (2004),
    Lecture notes on random walks in random environment.
    {\it \'Ecole d'\'et\'e de probabilit\'es de
    Saint-Flour 2001}. Lecture Notes in
    Math. {\bf 1837}, pp.~189--312. Springer, Berlin.

} 


\end{document}